\documentclass[a4paper,12pt]{article}
\usepackage{amsfonts,amssymb}
\title{ Topics on nonlinear generalized functions}

\author{J. F. Colombeau,\\
Institut Fourier, Universit\'e de Grenoble 1,   France\\ jf.colombeau@wanadoo.fr.}

\begin{document}
\maketitle

\begin{abstract} The aim of this paper is to give the text of a recent introduction to nonlinear generalized functions exposed in my talk in the congress GF2011, which was asked by several participants. Three representative topics were presented: two recalls "Nonlinear generalized functions  and their connections with distribution theory", "Examples of applications", and a recent  development: "Locally convex topologies and compactness: a functional analysis of nonlinear generalized functions".\\
\end{abstract}

\textit{ 1. Locally convex topologies and compactness}. \\

We start with this topic because of its complete novelty and its somewhat unexpected character which makes  the audience more interested in it than in the recalls. Let $\Omega$ be an open set in $\mathbb{R}^n $ and let $\mathcal{G}(\Omega)$ denote the special algebra of generalized functions on $\Omega$, i.e. there is no canonical inclusion of the vector space of distributions into $\mathcal{G}(\Omega)$: we can obtain such an inclusion after choice of a mollifier. A topology on $\mathcal{G}(\Omega)$ was defined in [2,1]. This topology is not a vector space topology.  Later this topology was rediscovered by D. Scarpalezos [16,17,13] who gave it the name of "sharp topology", and improved by Scarpalezos, Aragona, Garetto,Verwaerde, \dots.\\

Later I found a proof of nonexistence of a Hausdorff locally convex topology on $\mathcal{G}(\Omega)$ having some natural needed properties. At this point $\mathcal{G}(\Omega)$ permits concrete applications  such as in nonlinear elasticity and acoustics in heterogeneous media.  Since these methods are based on topology,  complete metric spaces (contractive mappings), or Hilbert spaces, or compactness,$\dots$ a topology is needed and one could strongly wish a classical topology so as to be able to use the well developped theory of locally convex spaces, Banach spaces, Hilbert spaces,$\dots$, as this has been done in the various spaces of the theory of distributions:  locally convex spaces in the Schwartz theory, Sobolev spaces,$\dots$. Therefore the natural character of the following question:\\

Does there exists a convenient subalgebra of $\mathcal{G}(\Omega)$ having the requested topological properties? We shall expose here that the answer is yes: even there are many of them.\\

\textit{Theorem.} There exists subalgebras of $\mathcal{G}(\Omega)$ with the following properties:\\
i) they are Hausdorff locally convex algebras in which all bounded sets are relatively compact,\\
ii) they contain "most" distributions on $\Omega$ and all bounded sets of distributions are relatively compact\\ 
iii) all partial derivatives are linear continuous from any such algebra into itself. Many nonlinearities are also continuous and internal.\\

The topological algebras we have constructed are nonmetrizable, complete, Schwartz and some of them nuclear. \\

The situation can be presented as follows: distribution theory provides a synthesis between differentiation and irregular functions: any irregular function is a distribution and any partial derivative of a distribution is still a distribution. $\mathcal{G}(\Omega)$ has added nonlinearities into this context: in particular any product of two elements of $\mathcal{G}(\Omega)$ is still an element of $\mathcal{G}(\Omega)$. The theory of distributions has also provided to mathematics a very rich variety of topologies having optimal properties (but limited to vector spaces due to the impossibility to multiply conveniently the distributions). Now the above subalgebras bring algebra topologies that are as rich as the topologies of the spaces of distributions, whith this basic difference that these topologies are compatible with nonlinearities, as this could be expected in the context of nonlinear generalized functions. In other word in these algebras one has 
$$ u_n\rightarrow u, v_n\rightarrow v \Rightarrow u_n.v_n \rightarrow u.v$$
exactly as in the classical spaces of $\mathcal{C}^\infty$ functions. Further in natural cases ii) means that 
$$ u_n\rightharpoonup u \Rightarrow u_n \rightarrow u$$
where the left one sided arrow means the classical (weak) convergence in spaces of distributions and the right complete arrow is the topological convergence in these algebras (there is no mistake in this implication as it could be expected at first sight).  Often one has a bounded family of approximate solutions of an equation; if one considers the strong topologies in some well chosen Banach spaces then the product $u_n v_n $ tends to $uv$ if $ u_n\rightarrow u$ and $ v_n\rightarrow v$, but from the Riesz theorem such bounded families are not a priori relatively compact so that one cannot extract convergent subsequences. To have the needed compactness tool one is usually forced to work with weak (or *weak) topologies. Then it is well known that these topologies are incompatible with a passage to the limit in nonlinear terms:
$$ u_n\rightharpoonup u , v_n\rightharpoonup v \not\Rightarrow u_n v_n \rightharpoonup uv$$

Therefore we can hope that this context could be useful in nonlinear problems.\\

As presented above the context announced in the theorem looks therefore like the Schwartz presentation of distributions in which the locally convex vector spaces would be replaced by locally convex topological algebras. Since the constructions of these topological algebras are completely different from the definition of the spaces in Schwartz theory there is work for those mathematicians who love topological vector spaces. This context can be explained  in the context  with the concept of Banach and Hilbert spaces only. For instance: \\ 

 There is an infinite  increasing sequence $(H_n)$ of separable Hilbert spaces with nuclear inclusion maps $H_n\rightarrow H_{n+1}$. The union of the $H_n$'s is strictly smaller than $\mathcal{G}(\Omega)$ but it contains significative spaces of distributions such that any bounded set in these spaces of distributions is contained and bounded in one $H_n$. One can choose these Hilbert spaces such that $H_n.H_n\subset H_{n+1}$ where $H_n.H_n:=\{x.y\}, x\in H_n, y\in H_n$, and such that $\frac{\partial}{\partial x_i}H_n\subset H_{n+1}.$\\
 
 Numerous questions were raised after the talk. The applications to PDEs have not yet been investigated. Other questions concerned the construction of these algebras: it is quite delicate and could not be given in full in such a short time. The idea of the construction I did is that in these subalgebras of $\mathcal{G}(\Omega)$, whose elements are of course  equivalence classes (from the definition of $\mathcal{G}(\Omega)$), there is in each equivalence class a privilieged representative so that there is an algebraic isomorphism between the algebra of equivalence classes and a classical algebra (i.e. without quotient) made of these privilieged representatives. The absence of a quotient in the algebra of the privilieged representatives permits to define there natural Hausdorff locally convex topologies that it suffices to transport to the subalgebra of $\mathcal{G}(\Omega)$ under concern. One uses deeply the theory of topological vector spaces [3, or other books] and the theory of nuclear spaces [12,15].\\

In conclusion it is likewise that the above opens the possibility of an original functional analysis of nonlinear generalized functions as rich as the classical linear functional analysis, but significantly different.\\

\textit{2. Nonlinear generalized functions as a natural  continuation of Schwartz presentation of distributions.} \\

The purpose of this section is to show to mathematicians who do not know the nonlinear generalized functions that their construction is indeed  naturally connected to the classical linear theory as presented by Schwartz [18]. For this we explain the construction of the nonlinear generalized functions for mathematicians in a pedagogical way to stress the basic facts and hide the (slightly) technical points that are minor but could be a drawback for a clear understanding of the situation. We adopt the Schwartz notation: $\mathcal{D}(\Omega)=\mathcal{C}_c^\infty(\Omega)=$ the space of all $\mathcal{C}^\infty$ functions on $\Omega$ with compact support, $\mathcal{D}'(\Omega)=$ the vector space of all distributions on $\Omega$, $\mathcal{E}'(\Omega)=$ the vector space of all distributions on $\Omega$ with compact support, $\delta_x=$ the Dirac measure centered at the point $x\in \Omega, \ \mathbb{K}= \mathbb{R}$ or $\mathbb{C}$.\\

In the algebra $\mathcal{C}^\infty(\mathcal{D}(\Omega))$ of all $\mathcal{C}^\infty$ maps $\Phi:\mathcal{D}(\Omega)\rightarrow \mathbb{K}$ one considers an equivalence relation $\mathcal{R}$, and -modulo some technical details- the algebra $\mathcal{G}(\Omega)$ is the set of all equivalence classes. To understand this equivalence relation -whose definition is somewhat technical- we describe it on subspaces of $\mathcal{C}^\infty(\mathcal{D}(\Omega))$  where it takes a particularly simple form.\\
$$\mathcal{C}^\infty(\mathcal{E}'(\Omega)) \subset \mathcal{C}^\infty(\mathcal{D}(\Omega))$$
through the restriction map $\Phi\rightarrow \Phi_{|\mathcal{D}(\Omega)}$ which is injective since $\mathcal{D}(\Omega)$ is dense in $\mathcal{E}'(\Omega)$. In $\mathcal{C}^\infty(\mathcal{E}'(\Omega))$ the equivalence relation $\mathcal{R}$ takes a particularly simple form:
$$\Phi_1\mathcal{R}\Phi_2 \Leftrightarrow \Phi_1(\delta_x)=\Phi_2(\delta_x) \forall x\in \Omega$$
In the subspace  $\mathcal{D}'(\Omega)$ of $\mathcal{C}^\infty(\mathcal{D}(\Omega))$ then the equivalence relation $\mathcal{R}$ reduces to the identity. That is why  $\mathcal{D}'(\Omega) \subset \mathcal{G}(\Omega)$  and why this equivalence relation  was not perceived by specialists of distribution theory. The statement of the equivalence relation in $\mathcal{C}^\infty(\mathcal{D}(\Omega))$ is merely an extension of the above in which $\delta_x$ is replaced by a net $\{\rho_{\epsilon,x}\}_\epsilon$ of elements of $\mathcal{D}(\Omega)$ that approximate $\delta_x$ in a standard way. On such a net the equality $\Phi(\delta_x)=0$ is replaced by a suitable fast decrease to 0 of the values $\Phi(\rho_{\epsilon,x})$ when $\epsilon\rightarrow 0$.\\

The Schwartz impossibility result states that "multiplication of distributions is impossible in any mathematical context possibly different from distribution theory". This appears to be in contradiction with the existence of the algebra $\mathcal{G}(\Omega)$. Here is a simplified  version of what appears to be  the "core" of this result. Let $H$ denote the Heaviside function.\\

$\int_{-\infty}^{+\infty}(H^2-H)H'dx=0$ because $H^2=H$. On the other hand\\
$$\int_{-\infty}^{+\infty}(H^2-H)H'dx=[\frac{H^3}{3}-\frac{H^2}{2}]_{-\infty}^{+\infty}=\frac{1}{3}-\frac{1}{2}=-\frac{1}{6}.$$

What is the correct result? If one admits that $H^2=H$ as adopted in classical mathematics from Lebesgue integration theory then the second line proves that multiplication of $H$ and $H'$ is impossible. In the theory of nonlinear generalized functions the second line is true, which shows that necessarily in this theory $H^2\not=H$. This is in perfect agreement with physical intuition: indeed $H$ can be viewed as an idealization of a continuous phenomenon that "jumps" on a very small interval, say of length $\epsilon$ around $x=0$. Then $H^2-H$ appears as having nonzero values   located on this small interval, and $H'$ appears to have large values located on this interval. When $\epsilon\rightarrow 0$ the first line of the above calculations has the form $0\times \infty$ whose undeterminacy here is solved by the second line. Therefore the classical algebra of step
 functions is not a faithful subalgebra of $\mathcal{G}(\mathbb{R})$. The Schwartz impossibility result extends this fact to the algebra $\mathcal{C}(\mathbb{R})$ of all continuous functions on $\mathbb{R}$. This fact requests explanations since now one has two different products of continuous functions, which is a priori unacceptable. Indeed the situation is  solved in a satisfactory way simply by noticing that the difference between the two products is always insignificant as long as one perfoms calculations that make sense within distribution theory. For instance it is true that $H^2$ and $H$ can be identified as long as one considers integration with a test function: $\forall \psi \in \mathcal{D}(\Omega)$ then $\int H \psi dx=\int H^2 \psi dx$. In $\mathcal{G}(\Omega)$ we say that $H^2$ and $H$ are associated, i.e. their integration on any test function in $\mathcal{D}(\Omega)$ gives same result. The association is another equivalence relation in $\mathcal{C}^\infty(\mathcal{D}(\Omega))$, less restrictive than the relation $\mathcal{R}$ considered above. For all calculations valid within distribution theory associated objects give same result and the restriction of the association to the space of distributions  reduces again to the equality of distributions. Therefore the fact noticed by Schwartz that he interpreted as impossibility to multiply the distributions does not cause problem: the calculations that make sense within the distributions give always the same result when reproduced in $\mathcal{G}(\Omega)$ This is developped in detail in numerous pedagogical texts  and talks such as [4-8,11,13,14]. \\

\textit{3. Examples of applications} This has already been developped in expository texts such as [6,9-11,14].\\


\begin{thebibliography}{<2>}

\bibitem{1bl12} H.A. Biagioni. A Nonlinear Theory of Generalized Functions. Lecture Notes in Math. 1421. Springer. 1990.  


\bibitem{1bl12} H.A. Biagioni, J.F. Colombeau. New Generalized Functions and $\mathcal{C}^\infty$ functions with values in generalized complex numbers. J. London Math. Soc. 2,33, 1986, p. 169-179.
\bibitem{1bl12} N. Bourbaki. Topological Vector Spaces.Hermann, Paris.

\bibitem{1bl12}J.F. Colombeau. New Generalized Functions and Multiplication of Distributions. North-Holland. 1984.
\bibitem{1bl12}J.F. Colombeau. Elementary Introduction to Nonlinear Generalized FunctionS. North-Holland. 1985.

\bibitem{1bl12} J.F. Colombeau. Multiplication of Distributions . Lecture Notes in Math. 1532. Springer. 1992.
\bibitem{1bl12} J.F. Colombeau. Generalized functions and infinitesimals. ArXiv 0610264.
\bibitem{1bl12} J.F. Colombeau. Generalized functions as a tool for nonsmooth nonlinear problems. ArXiv 061071.


\bibitem{1bl12} J.F. Colombeau, A Gsponer,B. Perrot. Generalized functions and the Heisenberg-Pauli formalism of QFT ArXiv 07052396


\bibitem{1bl12} J.F. Colombeau, A Gsponer. The Heisenberg-Pauli canonical formalism of QFT in the rigorous setting of nonlinear generalized functions ArXiv 07083425.

\bibitem{1bl12} M. Grosser, M. Kunzinger, M. Oberguggenberger, R. Steinbauer. Geometric Theory of Generalized Functions with Applications to General Relativity. Kluwer 2001.
\bibitem{1bl12} A. Grothendieck. Produits tensoriels topologiques; Memoirs of the AMS.
\bibitem{1bl12} M. Nedeljkov, S. Pilipovic, D. Scarpalezos. The Linear Theory of Colombeau Generalized Functions. Pitman Research Notes in Math. 1998.
\bibitem{1bl12} M. Oberguggenberger. Multiplication of Distributions and Applications to Partial Differential Equations. Pitman Research Notes in Math. 1992.
\bibitem{1bl12} Pietsch. Nuclear locally convex  spaces. Erg. der Math. 66 , Springer Verlag, 1972.
\bibitem{1bl12} D. Scarpalezos.  Colombeau generalized functions, topological structures, microlocal properties. A simplified point of view.
\bibitem{1bl12} D. Scarpalezos. Private communication, about 1989.
\bibitem{1bl12} L. Schwartz. Theorie des distributions. Hermann, Paris, numerous editions.

\end{thebibliography}
\end{document}